\documentstyle[10pt]{article}
\begin{document}
\begin{center}
{\Large \bf  Krull dimension of tensor products of pullbacks
}\vspace{1mm}
\end{center}
\begin{center}
{\bf Samir Bouchiba}\\
{\small \it Department of Mathematics, University Moulay Ismail, Meknes 50000, Morocco} \vspace{8mm}
\end{center}

\noindent {\small {\bf Abstract}}

{\small This paper is concerned with the study of the dimension theory of tensor products of
algebras over a field $k$. We answer an open problem set
in [6] and
compute dim$(A\otimes_kB)$ when $A$ is a $k$-algebra arising from a
specific pullback construction involving AF-domains and $B$ is an arbitrary $k$-algebra. On the other hand, we deal with
the question (Q) set in [5] and show, in particular, that
such a pullback $A$ is in fact a generalized AF-domain.}
\vspace{5mm}\bigskip

\noindent {\bf 1. Introduction}\\

All rings considered in this paper are commutative with identity
element and all ring homomorphisms are unital. Throughout, $k$
stands for a field. We shall use t.d.($A:k)$, or t.d.($A)$ when no
confusion is likely, to denote the transcendence degree of a
$k$-algebra $A$ over $k$, $A[n]$ to denote the polynomial ring
$A[X_1,...,X_n]$ and $p[n]$ to denote the prime ideal
$p[X_1,...,X_n]$ of $A[X_1,...,X_n]$ for each prime ideal $p$ of
$A$. Also, we use Spec($A)$ to denote the set of prime ideals of a
ring $A$ and $\subset$ to denote proper set inclusion. All
$k$-algebras considered throughout this paper are assumed to be of
finite transcendence degree over $k$. Any unreferenced material is
standard as in [11], [15], [16] and [17].

Several authors have been interested in studying the prime ideal
structure and related topics of tensor products of algebras over a
field $k$. The initial impetus for these investigations was a paper
of R. Sharp on Krull dimension of tensor products of two extension
fields. In fact, in [19], Sharp proved that, for any two extension
fields $K$ and $L$ of $k$, dim$(K\otimes_kL)=$ min(t.d.($K),$
t.d.($L))$ (actually, this result appeared ten years earlier in
Grothendieck's EGA [13, Remarque 4.2.1.4, p. 349]). This formula is
rather surprising since, as one may expect, the structure of the
tensor product should reflect the way the two components interact
and not only the structure of each component. This fact is what most
motivated Wadsworth's work in [20] on this subject. His aim was to
seek geometric properties of primes of $A\otimes_kB$ and to widen
the scope of algebras $A$ and $B$ for which dim$(A\otimes_kB)$
depends only on individual characteristics of $A$ and $B$. The
algebras which proved tractable for Krull dimension computations
turned out to be those domains $A$ which satisfy the altitude
formula over
$k$ (AF-domains for short), that is,\\
\noindent {\large \bf ---------------}

{\footnotesize {\it Mathematics Subject Classification} (2000):
Primary 13C15; Secondary 13B24.

{\it E-mail address}: sbouchiba@hotmail.com

{\it keywords:} Krull dimension, tensor product, prime ideal,
AF-domain.}\newpage

 $$ht(p)+\mbox {t.d.}(\frac Ap)=\mbox { t.d.}(A)$$ for all
prime ideals $p$ of $A$. It is worth noting that the class of
AF-domains contains the most basic rings of algebraic geometry,
including finitely generated $k$-algebras that are domains.
Wadsworth proved, via [20, Theorem 3.8], that if $A_1$ and $A_2$ are
AF-domains, then
$$\mbox {dim}(A_1\otimes_kA_2)=\mbox { min}\Big {(}\mbox
{dim}(A_1)+\mbox {t.d.}(A_2),\mbox { t.d.}(A_1)+\mbox {dim}(A_2)\Big
{)}.$$ His main theorem stated a formula for dim($A\otimes_kB)$
which holds for an AF-domain $A$, with no restriction on $B$,
namely:$$ \begin{array}{lll} \mbox {dim}(A\otimes_kB)&=&D\Big
{(}\mbox {t.d.}(A),\mbox {dim}(A),B\Big
{)}\\
&:=&\mbox {max}\Big {\{}ht(q[\mbox {t.d.}(A)])+\mbox {min}\Big
{(}\mbox {t.d.}(A),\mbox { dim}(A)+\displaystyle {\mbox {t.d.}(\frac
Bq})\Big {)}:q\in \mbox {Spec}(B)\Big {\}}\\
&&\mbox {[20, Theorem 3.7].}
\end{array}$$

Our aim in [6] was to extend Wadsworth's results in a different way,
namely to tensor products of $k$-algebras arising from pullbacks. In
this regard, we use previous deep investigations on prime ideal
structure of various pullbacks, as in [1]. Our main result in [6]
states the following:\\

{\it Let $T_i$ be a $k$-algebra which is an integral domain and
$M_i$ a maximal ideal of $T_i$ such that $ht(M_i)=$ dim$(T_i)$,
$K_i=\displaystyle {\frac {T_i}{M_i}}$, $\varphi_i$ the canonical
surjection from $T_i$ to $K_i$, $D_i$ a subring of $K_i$, and
$A_i=\varphi_i^{-1}(D_i)$ be the issued pullback, for $i=1,2$.
Assume that $T_i$ and $D_i$ are AF-domains, for $i=1,2$. Then
$$\mbox {dim}(A_1\otimes_kA_2)=\mbox { max}\Big {\{}ht(M_1[\mbox
{t.d.}(A_2)])+D\Big {(}\mbox {t.d.}(D_1),\mbox { dim}(D_1),R_2\Big
{)},$$ $$ht(M_2[\mbox {t.d.}(A_1)])+D\Big {(}\mbox {t.d.}(D_2),\mbox
{ dim}(D_2),R_1\Big {)}\Big {\}},$$} where $D(s,d,A):=$ max$\Big
{\{}ht(p[s])+$min$\Big {(}s,d+$t.d.$(\displaystyle {\frac Ap})\Big
{)}:p\in$ Spec$(A)\Big {\}}$ for any $k$-algebra $A$ and any
positive integers $0\leq d\leq s$. This theorem allows one to
compute the Krull dimension of tensor product of two $k$-algebras
for a large family of (not necessarily AF-domains) $k$-algebras.
Further, we set in [6] the open problem of computing
dim$(A_1\otimes_kA_2)$ when only $T_1$ and $D_1$ are assumed to be
AF-domains.

On the other hand, in [14], Jaffard proved that, for any ring $A$
and any positive integer $n$, the Krull dimension of $A[n]$ can be
realized as the length of a special chain of $A[n]$. Recall that a
chain $C=\{Q_0\subset Q_1\subset ...\subset Q_s\}$ of prime ideals
of $A[n]$ is called a special chain if for each $Q_i$, the ideal
$(Q_i\cap A)[n]$ belongs to $C$. Subsequently, Brewer et al. gave an
equivalent and simple version of Jaffard's theorem. Actually, they
showed that, for each positive integer $n$ and each prime ideal $P$
of $A[n]$, $ht(P)=ht(q[n])+ht(\displaystyle {\frac P{q[n]})}$ [10,
Theorem 1], where $q:=P\cap A$. Taking into account the natural
isomorphism $B[n]\cong k[n]\otimes_kB$ for each $k$-algebra $B$, we
generalized in [6] this special chain theorem to tensor products of
$k$-algebras. Effectively, we proved that if $A$ and $B$ are
$k$-algebras such that $A$ is an AF-domain, then for each prime
ideal $P$ of $A\otimes_kB$,
$$ht(P)=ht(A\otimes_kq)+ht(\frac P{A\otimes_kq})=ht(q[\mbox
{t.d.}(A)])+ht(\frac P{A\otimes_kq}),$$ where $q=P\cap B$  (cf. [6,
Lemma 1.5]). It turned out that this very geometrical property
totally characterizes the AF-domains. In fact, we proved, in [4],
that the following statements are equivalent for a domain $A$ which
is a $k$-algebra:\\

{\it a) $A$ is an AF-domain;

b) $A$ satisfies SCT (for special chain theorem), that is, for each
$k$-algebra $B$ and each prime ideal $P$ of $A\otimes_kB$ with
$q:=P\cap B$,
$$ht(P)=ht(q[\mbox {t.d.}(A)])+ht(\displaystyle {\frac
P{A\otimes_kq})}=ht(A\otimes_kq)+ht(\frac P{A\otimes_kq}) \mbox {
[4, Theorem 1.1]}.$$}In view of this, we generalized in [5] the
AF-domain notion by setting the following definitions:

{\it We say that a $k$-algebra $A$ satisfies GSCT (for generalized
special chain theorem) with respect to a $k$-algebra $B$ if
\begin{center}
$ht(P)=ht(p\otimes_kB+A\otimes_kq)+ht(\displaystyle {\frac
P{p\otimes_kB+A\otimes_kq})}$\\

 for each prime ideal $P$ of
$A\otimes_kB$, with $p=P\cap A$ and $q=P\cap B$,
\end{center}

\noindent and we call a generalized AF-domain (GAF-domain for short)
a domain $A$ such that $A$ satisfies GSCT with respect to any
$k$-algebra $B$.}\\ There is no known example in the literature of a
$k$-algebra $A$ which is a domain and which is not a GAF-domain.
This may lead one to ask whether any $k$-algebra which is a domain
is a GAF-domain. We were concerned in [5] with the following question:\\

(Q): Is any domain $A$ which is a $k$-algebra such that the
polynomial ring $A[n]$ is an AF-domain, for some positive integer
$n$, a GAF-domain?\\

\noindent We gave in [5] partial results settling in the affirmative
the above question (Q). First, we proved that an AF-domain $A$ is in
fact a GAF-domain, thus in particular, any finitely generated
algebra over $k$ which is a domain is a GAF-domain. Also, we proved
that (Q) has a positive answer in the case where $A$ is
one-dimensional. Our main result in [5] tackles the case
$n=1$ of $(Q)$. It computes dim$(A\otimes_kB)$ for a $k$-algebra $A$
such that $A[X]$ is an AF-domain and for an arbitrary $k$-algebra
$B$ generalizing Wadsworth's main theorem [20, Theorem 3.7] and
further asserts that $A$ is a GAF-domain. We ended that paper by an
example of a GAF-domain $A$ such that, for any positive integer $n$,
the polynomial ring $A[n]$ is not an AF-domain.

Our objective in this paper is twofold. On the one hand, we handle
the above-mentioned problem set in [6] and compute
dim$(A\otimes_kB)$ when $A$ is a pullback arising from the above
construction and $B$ is an arbitrary $k$-algebra. On the other hand,
we prove that the answer to the question (Q) set in [5] is
affirmative for such a pullback construction $A$. Besides, our main
result, Theorem 2.8, is, in particular, an important step towards
determining a general formula for dim$(A\otimes_kB)$ in the case
where $A[n]$ is an AF-domain for some positive integer $n$ and $B$
is an arbitrary $k$-algebra. It states the following:\\

{\it Let $T$ be a $k$-algebra which is a domain and $M$ a maximal
ideal of $T$. Let $K=\displaystyle {\frac TM}$ and $D$ be a subring
of $K$. Let $\varphi:T\longrightarrow K$ be the canonical surjective
homomorphism and $A:=\varphi^{-1}(D)$. Assume that $T$ and $D$ are
AF-domains and $T_M$ is catenarian. Then, $A$ is a GAF-domain and
for an arbitrary $k$-algebra $B$,
$$\mbox {dim}(A\otimes_kB)=\mbox { max}\Big {\{}D\Big (\mbox
{t.d.}(A),d,B\Big );ht(M)+\mbox {max}\Big {\{}ht(q_1[\mbox
{t.d.}(A)])+\displaystyle { ht\Big (\frac q{q_1}[\mbox
{t.d.}(D)]\Big )}+$$
$$\displaystyle {\mbox {min\Big {(}t.d.}(\frac B{q_1}),\mbox {
t.d.}(K:D)\Big {)}}+\displaystyle {\mbox {min\Big {(}t.d.}(D),\mbox
{ dim}(D)+\mbox {t.d.}(\frac Bq)\Big {)}}:q_1\subseteq q\in\mbox {
Spec}(B)\Big {\}}\Big {\}},$$ where $d:=\mbox { sup}\Big
{\{}ht(Q):Q\in$ Spec$(T)$ with $M\not\subset Q\Big {\}}$}.\\

\noindent Direct consequences of this theorem are provided as well
as a case where we may drop the catenarity property of $T_M$ is
exhibited. An example to illustrate our findings closes this paper.

Recent developments on height and grade of (prime) ideals as well as
on dimension theory in tensor products of $k$-algebras are to be
found in [2-7]. Concerning the study of the transfer to tensor products of algebras
 of the S-property, strong S-property, and catenarity,
we refer the reader to [8].\\

\noindent {\bf 2. Main results}\\

First, for the convenience of the reader, we catalog some basic
facts and results connected with the tensor product of
$k$-algebras. These will be used frequently in the sequel without explicit mention.\\

Let $A$ and $B$ be two $k$-algebras. If $p$ is a prime ideal of $A$,
$r=$ t.d.($\displaystyle {\frac Ap})$ and $\overline
{x_1},...,\overline {x_r}$ are elements of $\displaystyle {\frac
Ap}$, algebraically independent over $k$, with the $x_i \in A$, then
it is easily seen that $x_1,...,x_r$ are algebraically independent
over $k$ and $p\cap S=\emptyset$, where
$S=k[x_1,...,x_r]\setminus\{0\}$. If $A$ is an integral domain, then
$ht(p)+$t.d.($\displaystyle {\frac Ap})\leq$ t.d.($A)$ for each
prime ideal $p$ of $A$ (cf. [21, p. 37] ). Now, assume that $S_1$
and $S_2$ are multiplicative subsets of $A$ and $B$, respectively,
then $S_1^{-1}A\otimes _kS_2^{-1}B\cong S^{-1}(A\otimes _kB)$, where
$S=\{s_1\otimes s_2:s_1\in S_1$ and $s_2\in S_2\}$. We assume
familiarity with the natural isomorphisms for tensor products. In
particular, we identify $A$ and $B$ with their respective images in
$A\otimes_kB$. Also, $A\otimes_kB$ is a free (hence faithfully flat)
extension of $A$ and $B$. Moreover, recall that an AF-domain $A$ is
a locally Jaffard domain, that is, $ht(p[n])=ht(p)$ for each prime
ideal $p$ and each positive integer $n$ [20, Corollary 3.2].
Finally, we refer the reader to the useful result of Wadsworth [20,
Proposition 2.3] which yields a classification of the prime
ideals of $A\otimes_kB$ according to their contractions to $A$ and $B$.\\

We begin by recalling from [3], [5], [6] and [20] the following useful results.\\

\noindent {\bf Proposition 2.1 [6, Lemma 1.3].} {\it Let $A$ and $B$
be $k$-algebras such that $B$ is a domain. Let $p$ be a prime ideal
of $A$. Then, for each prime ideal $P$ of $A\otimes_kB$ which is
minimal over $p\otimes_kB$,
$$ht(P)=ht(p\otimes_kB)=ht(p[\mbox {t.d.}(B)]).$$}

\noindent {\bf Proposition 2.2 [20, Proposition 2.3].} {\it Let $A$
and $B$ be $k$-algebras and let $p\subseteq p^{\prime}$ be prime
ideals of $A$ and $q\subseteq q^{\prime}$ be prime ideals of $B$.
Then the natural ring homomorphism $\varphi: \displaystyle {\frac
{A\otimes_kB}{p\otimes_kB+A\otimes_kq}\longrightarrow \frac
Ap\otimes_k\frac Bq}$ such that $\varphi (\overline
{a\otimes_kb})=\overline a\otimes_k\overline b$ for each $a\in A$
and each $b\in B$, is an isomorphism and $$\varphi (\displaystyle
{\frac
{p^{\prime}\otimes_kB+A\otimes_kq^{\prime}}{p\otimes_kB+A\otimes_kq})=\frac
{p^{\prime}}p\otimes_k\frac Bq+\frac Ap\otimes_k\frac {q^{\prime}}q}.$$}\\

\noindent {\bf Proposition 2.3 [6, Lemma 1.2 and Proposition 2.2].}
{\it Let $T$ be an integral domain which is a $k$-algebra, $M$ a
maximal ideal of $T$, $K:=\displaystyle {\frac TM}$ and $\varphi:
T\rightarrow K$ the canonical surjective homomorphism. Let $D$ be a
proper subring of $K$ and $A:=\varphi^{-1}(D)$. Assume that $T$ and
$D$ are AF-domains. Then

1) The polynomial ring $A[$t.d.$(K:D)]$ is an AF-domain.

2) $ht(p[n])=ht(p)+$min$\Big {(}n,$ t.d.$(K:D)\Big {)}$ for each prime ideal $p$ of
$A$ such that $M\subseteq p$ and each positive integer $n$.}\\

The following easy result is probably well known. We refer the  reader to [5] for a detailed proof.\\

\noindent {\bf Proposition 2.4 [5, Proposition 2.4].} {\it Let $A$ be a ring. Let
$I\subseteq J$ be ideals in $A$. Then $$ht(I)+ht(\displaystyle
{\frac JI})\leq ht(J).$$}

\noindent {\bf Proposition 2.5.} {\it Let $A$ and $B$ be
$k$-algebras such that $A$ is an AF-domain. Let $P\in$
Spec$(A\otimes_kB)$, $p=P\cap A$ and $q=P\cap B$. Then

1) $ht(P)=ht(\displaystyle { A\otimes_kq)+ht(\frac P{A\otimes_kq})}$
[6, Lemma 1.5].

2) $ht(P)=ht(q[\mbox {t.d.}(A)])+ht(p)+ht(\displaystyle {\frac
P{p\otimes_kB+A\otimes_kq}}).$}

\noindent {\bf Proof.} 2) As $\displaystyle {\frac
P{A\otimes_kq}\cap \frac Bq=(\overline 0)}$, $\displaystyle
{\frac P{A\otimes_kq}}$ survives in $A\otimes_kk_B(q)$, where
$k_B(q)$ denotes the quotient field of $\displaystyle {\frac Bq}$,
so that, applying (1), we get\\

$\begin{array}{lll} ht(\displaystyle {\frac
P{A\otimes_kq}})&=&ht\Big (p\Big [$t.d.$\Big (\displaystyle {\frac Bq)\Big ]\Big )+ht\Big (\frac
{P/(A\otimes_kq)}{p\otimes_k\frac Bq}\Big )}\\
&=&ht(p)+ht(\displaystyle {\frac P{p\otimes_kB+A\otimes_kq}}),$ as
$A$ is a locally Jaffard domain$\\
&&\\
&&$and $p\otimes_k\displaystyle {\frac Bq\cong \frac {p\otimes_kB+A\otimes_kq}{A\otimes_kq}}$, by Proposition 2.2$.
\end{array}
$\\

\noindent Hence $ht(P)=ht(q[$t.d.$(A)])+ht(p)+ht(\displaystyle
{\frac
P{p\otimes_kB+A\otimes_kq}})$, as desired. $\Box$\\

A domain $A$ is said to be catenarian if for each chain of prime
ideals $p\subseteq q$ of $A$, $$ht(p)+ht(\displaystyle {\frac
qp})=ht(q) \mbox { (cf. [9])}.$$

\noindent {\bf Proposition 2.6.} {\it Let $A$ be an AF-domain. If
$A$ is catenarian, then $\displaystyle {\frac Ap}$ is an AF-domain
for each
prime ideal $p$ of $A$.}\\

\noindent {\bf Proof.} Assume that $A$ is catenarian and fix $p\in$
Spec$(A)$. Let $p\subseteq q\in$ Spec$(A)$. Then\\

$\begin{array}{lll} ht(\displaystyle {\frac qp)+\mbox {t.d.}(\frac
Aq)}&=&\displaystyle {ht(q)-ht(p)+\mbox {t.d.}(\frac Aq)}\\
&=&\mbox{t.d.}(A)-ht(p)\mbox { as }A\mbox { is an AF-domain}\\
&=&\displaystyle {\mbox {t.d.}(\frac Ap)}.
\end{array}
$\\

\noindent Hence $\displaystyle {\frac Ap}$ is an AF-domain. $\Box$\\

 Let $A$ and $B$ be $k$-algebras and $P$ be a prime ideal of
$A\otimes _kB$. Let
 $q_0\in$ Spec$(B)$ such that $q_0\subset P\cap B$. We denote by
$\lambda \Big {(}(.,q_0),P\Big {)}$ the maximum of lengths of chains
of prime ideals of $A\otimes_kB$ of the form $P_0\subset P_1\subset
...\subset P_s=P$ such that $P_i\cap B=q_0$, for $i=0,1,...,s-1$.\\

\noindent {\bf Proposition 2.7 [3, Lemma 2.4].} {\it  Let $A$ and
$B$ be $k$-algebras and $P$ be a prime ideal of $A\otimes _kB$ with $p=P\cap A$ and $q=P\cap B$.
Assume that $A$
 and $B$ are integral domains. Then} $$\lambda {\Big (}(.,(0)),P{\Big )}\leq\mbox { t.d.}(A)-\mbox
{t.d.}(\frac Ap)+ht\Big (q\Big [\mbox {t.d.}(\frac Ap)\Big ]\Big )+ht(\frac
P{p\otimes_kB+A\otimes_kq}).$$

Finally, recall that, if $A$ is a $k$-algebra and $n\geq 0$ is an integer,
then the polynomial ring $A[n]$ is an AF-domain if and only if, for each prime ideal $p$ of $A$,
$$ht(p[n])+\mbox {t.d.}\displaystyle {(\frac Ap)}=\mbox {
t.d.}(A)\mbox { [6, Lemma 2.1]}.$$

Next, we announce the main theorem of this paper. It gives an answer
to the above-mentioned problem set in [6] as well as to the question
(Q) of [5] and represents an important step towards determining a
formula for dim$(A\otimes_kB)$ when $A[n]$ is an AF-domain for some
positive integer $n$, and
$B$ is an arbitrary $k$-algebra.\\

\noindent {\bf Theorem 2.8.} {\it Let $T$ be a $k$-algebra which is
a domain and $M$ a maximal ideal of $T$. Let $K=\displaystyle {\frac
TM}$ and $D$ be a subring of $K$. Let
$\varphi:T\longrightarrow K$ be the canonical surjective
homomorphism and $A:=\varphi^{-1}(D)$. Assume that $T$ and $D$ are
AF-domains and $T_M$ is
catenarian. Let $B$ be an arbitrary $k$-algebra and let $P\in$ Spec$(A\otimes_kB)$, $p=P\cap A$ and $q=P\cap B$. Then the following statements hold:\\

1) If $M\not\subset p$, then
$$ht(P)=\mbox ht(p)+ht(q[\mbox {t.d.}(A)])+\displaystyle {ht(\frac P{p\otimes_kB+A\otimes_kq})}.$$

2) If $M\subseteq p$, then $$ht(P)=ht(p)+\mbox {max}\Big
{\{}ht(q_1[\mbox {t.d.}(A)])+\displaystyle {ht\Big (\frac
q{q_1}[\mbox {t.d.}(D)]\Big )}+$$
$$\displaystyle {\mbox {min}\Big {(}\mbox {t.d.}(\frac B{q_1}),\mbox
{ t.d.}(K:D)\Big {)}}:\displaystyle {q_1\subseteq q\in \mbox {
Spec}(B)\Big {\}}+ht(\frac P{p\otimes_kB+A\otimes_kq})}.$$

3) $\mbox {dim}(A\otimes_kB)=\mbox { max}\Big {\{}D\Big (\mbox
{t.d.}(A),d,B\Big );ht(M)+\mbox {max}\Big {\{}ht(q_1[\mbox
{t.d.}(A)])+\displaystyle { ht\Big (\frac q{q_1}[\mbox {t.d.}(D)]\Big )}$

\hspace {1cm}$\displaystyle {\mbox {min\Big {(}t.d.}(\frac
B{q_1}),\mbox { t.d.}(K:D)\Big {)}}+\displaystyle {\mbox {min\Big
{(}t.d.}(D),\mbox { dim}(D)+\mbox {t.d.}(\frac Bq)\Big {)}}:q_1\subseteq q\in\mbox { Spec}(B)\Big {\}}\Big {\}},$\\
where $d:=\mbox { sup}\Big {\{}ht(Q):Q\in$ Spec$(T)$ with $M\not\subset
Q\Big {\}}$.\\

4) $A$ is a GAF-domain.}\\

\noindent {\bf Proof.} 1) Let $M\not\subset p$. Then, by [1, Lemma 2.1], there exists $p^{\prime}\in$ Spec$(T)$
 such that $p^{\prime}\cap A=p$, and $p^{\prime}$ satisfies $A_p=T_{p^{\prime}}$. Thus $A_p$ is
an AF-domain, so that by Proposition 2.5, $$ht(P)=\displaystyle
{ht(q[\mbox {t.d.}(A)])+ht(p)+ht(\frac
P{p\otimes_kB+A\otimes_kq})},\mbox { as desired}.$$

2) Assume that $M\subseteq p$. Then t.d.$\displaystyle {(\frac
Ap)=\mbox {t.d.}\Big (\frac {A/M}{p/M}\Big )\leq\mbox { t.d.}(D)}$, and
applying [1, Lemma 2.1], there exists $Q\in$ Spec$(D)$ such that
$p=\varphi^{-1}(Q)$ and the following diagram
$$\begin{array}{lll}
A_p&\longrightarrow& D_Q\\
&&\\
\downarrow&&\downarrow\\
&&\\
T_M&\longrightarrow&K
\end{array}
$$
is a pullback diagram. Therefore, as $ht(P)=ht(P(A_p\otimes_kB))$, we may
assume without loss of generality that $(T,M)$ is a quasilocal
catenarian domain, and thus $M$ is a divided prime ideal of $A$, via
[1, Lemma 2.1].\\
First, note that, by Proposition 2.1 and Proposition 2.2, we have,
$\forall q_1\subseteq
q\in$ Spec$(B)$,\\

$\left \{
\begin{array}{lll} ht(q_1[\mbox
{t.d.}(A)])&=&ht(A\otimes_kq_1)\\
&&\\
ht\Big (p\Big [$t.d.$(\displaystyle {\frac B{q_1}})\Big ]\Big )&=&ht(\displaystyle
{\frac {p\otimes_kB+A\otimes_kq_1}{A\otimes_kq_1}})\\
&&\\
\displaystyle {ht\Big (\frac q{q_1}[\mbox {t.d.}(D)]\Big )}&=&\displaystyle
{ht(\frac {M\otimes_kB+A\otimes_kq}{M\otimes_kB+A\otimes_kq_1})}.
\end{array}
\right.$\\

\noindent It follows, by Proposition 2.1, Proposition 2.2, Proposition 2.3 and
Proposition 2.4, that, $\forall q_1\subseteq q\in$ Spec$(B)$,
$$ht(p)+ht(q_1[\mbox {t.d.}(A)])+\displaystyle {ht\Big (\frac q{q_1}[\mbox {t.d.}(D)]\Big )}+
\mbox {min}\Big {(}\mbox {t.d.}(\displaystyle {\frac B{q_1}),\mbox {
t.d.}(K:D)}\Big {)}=$$
$$ht(M)+\mbox {min}\Big {(}\mbox {t.d.}(\displaystyle {\frac B{q_1}),\mbox {
t.d.}(K:D)}\Big {)}+ht(\displaystyle {\frac
pM})+ht(A\otimes_kq_1)+\displaystyle {ht(\frac
{M\otimes_kB+A\otimes_kq}{M\otimes_kB+A\otimes_kq_1})}=$$
$$ht\Big (M\Big [\mbox {t.d.}(\frac
B{q_1})\Big ]\Big )+ht(A\otimes_kq_1)+ht\Big (\displaystyle {\frac pM}\Big [\mbox
{t.d.}(\frac B{q})\Big ]\Big )+\displaystyle {ht(\frac
{M\otimes_kB+A\otimes_kq}{M\otimes_kB+A\otimes_kq_1})}=$$
$$\displaystyle {ht(\frac
{M\otimes_kB+A\otimes_kq_1}{A\otimes_kq_1})+ht(A\otimes_kq_1)+ht(\frac
{p\otimes_kB+A\otimes_kq}{M\otimes_kB+A\otimes_kq})}+\displaystyle
{ht(\frac
{M\otimes_kB+A\otimes_kq}{M\otimes_kB+A\otimes_kq_1})}\leq$$
$$ht(p\otimes_kB+A\otimes_kq).$$ Therefore, for each prime ideal $q_1\subseteq q$ of
$B$, $$ht(p)+ht(q_1[\mbox { t.d.}(A)])+\displaystyle {ht\Big (\frac
q{q_1}[\mbox {t.d.}(D)]\Big )}+\displaystyle {\mbox {min}\Big {(}\mbox
{t.d.}(\frac B{q_1}),\mbox { t.d.}(K:D)\Big {)}}+$$ $$\displaystyle
{ht(\frac P{p\otimes_kB+A\otimes_kq})}\leq ht(p\otimes_kB+A\otimes_kq)+ht(\displaystyle {\frac
P{p\otimes_kB+A\otimes_kq}})\leq ht(P),$$ establishing the
direct inequality. The proof of the reverse inequality falls into
the following two steps.\\

\noindent {\bf  \underline {Step 1}.} $B$ is an integral domain.\\
Our argument uses induction on dim$(T)$,
$ht(p)$ and $ht(q)$. First, note that
$$(*)\hspace {1cm} \mbox {max}\Big {\{}ht(q[\mbox
{t.d.}(A)])+ht\Big (p\Big [\displaystyle {\mbox {t.d.}(\frac Bq)\Big ]\Big ),ht(p[\mbox
{t.d.}(B)])+ ht\Big (q\Big [\mbox {t.d.}(\frac Ap)\Big ]\Big )\Big {\}}}\leq$$
$$ht(p)+\mbox {max}\Big {\{}ht(q_1[\mbox { t.d.}(A)])+\displaystyle {ht\Big (\frac q{q_1}[(\mbox {t.d.}(D)]
\Big )}+$$ $$\displaystyle {\mbox
{min}\Big {(}\mbox {t.d.}(\frac B{q_1}),\mbox { t.d.}(K:D)\Big
{)}}:q_1\subseteq q\in \mbox { Spec}(B)\Big {\}},$$ it suffices to
take $q_1=q$ and $q_1=(0)$. If either dim$(T)=0$ or $ht(p)=0$, then $T=K$ is a field and, thus, $A=D$
is an AF-domain, and applying Proposition 2.5 and $(*)$, we obtain the formula. Also, the case $ht(q)=0$ is fairly
easy via Proposition 2.5 and $(*)$. Then, assume that
dim$(T)>0$, $ht(p)>0$ and $ht(q)>0$. Consider a chain of prime
ideals $Q_0\subset Q_1\subset...\subset Q_h=P$ of $A\otimes_kB$ such
that $h=ht(P)$. Let $r:=$ max$\{m:Q_m\cap A\subset p$ or $Q_m\cap
B\subset q\}$. Let $Q=Q_r$, $p^{\prime}=Q_r\cap A$ and
$q^{\prime}=Q_r\cap B$. Hence, $ht(P)=ht(Q)+ht(\displaystyle {\frac
PQ})$. We are led to
discuss the following cases.\\

\noindent {\bf Case 1.}  $p^{\prime}=p$. Then, $q^{\prime}\subset
q$, and by inductive assumptions,
$$ht(Q)=ht(p)+\mbox {max}\Big
{\{}ht(q_1[\mbox { t.d.}(A)])+\displaystyle {ht\Big {(}\frac
{q^{\prime}}{q_1}[\mbox {t.d.}(D)]\Big )}+$$
$$\displaystyle {\mbox {min}\Big {(}\mbox {t.d.}(\frac B{q_1}),\mbox
{ t.d.}(K:D)\Big {)}}:\displaystyle {q_1\subseteq q^{\prime}\in
\mbox { Spec}(B)\Big {\}}}+ht(\displaystyle {\frac
Q{p\otimes_kB+A\otimes_kq^{\prime}}}),\mbox { and thus}$$
$$ht(P)\leq ht(p)+\mbox {max}\Big
{\{}ht(q_1[\mbox { t.d.}(A)])+\displaystyle {ht\Big (\frac
{q^{\prime}}{q_1}[\mbox {t.d.}(D)]\Big )}+$$
$$\displaystyle {\mbox {min}\Big {(}\mbox {t.d.}(\frac B{q_1}),\mbox
{ t.d.}(K:D)\Big {)}}:\displaystyle {q_1\subseteq q^{\prime}\in
\mbox { Spec}(B)\Big {\}}}+ht(\displaystyle {\frac
P{p\otimes_kB+A\otimes_kq^{\prime}}}).$$ As $\displaystyle {\frac
P{p\otimes_kB+A\otimes_kq^{\prime}}\cap \frac Ap}=(\overline 0)$, we
get by Proposition 2.5,
$$\begin{array}{lll}
ht(\displaystyle {\frac
P{p\otimes_kB+A\otimes_kq^{\prime}})}&=&\displaystyle {ht\Big (\frac
q{q^{\prime}}\Big [\mbox {t.d.}(\frac Ap)\Big ]\Big )+ht\Big {(}\frac
{P/(p\otimes_kB+A\otimes_kq^{\prime})}{(A/p)\otimes_k(q/q^{\prime})}\Big {)}}\\
&&\\
&=&\displaystyle {ht\Big (\frac q{q^{\prime}}\Big [\mbox {t.d.}(\frac
Ap)\Big ]\Big )+ht(\frac P{p\otimes_kB+A\otimes_kq})},
\end{array}$$
since, by Proposition 2.2, $\displaystyle {\frac Ap\otimes_k\frac
q{q^{\prime}}\cong \frac
{p\otimes_kB+A\otimes_kq}{p\otimes_kB+A\otimes_kq^{\prime}}}$. It
follows that\\

$\begin{array}{lll} ht(P)&\leq&ht(p)+\mbox {max}\Big
{\{}ht(q_1[\mbox { t.d.}(A)])+\displaystyle {ht\Big (\frac
{q^{\prime}}{q_1}[\mbox {t.d.}(D)]\Big )}+\\
&&\displaystyle {\mbox {min}\Big {(}\mbox {t.d.}(\frac B{q_1}),\mbox
{ t.d.}(K:D)\Big {)}}:\displaystyle {q_1\subseteq q^{\prime}\in
\mbox { Spec}(B)\Big {\}}}+\displaystyle {ht\Big (\frac
q{q^{\prime}}\Big [\mbox {t.d.}(\frac Ap)\Big ]\Big )+}\\
&&\displaystyle {ht(\frac
P{p\otimes_kB+A\otimes_kq})}\\
&&\\
&\leq&ht(p)+\mbox {max}\Big {\{}ht(q_1[\mbox {
t.d.}(A)])+\displaystyle {ht\Big (\frac {q^{\prime}}{q_1}[\mbox {t.d.}(D)]\Big )}+
\displaystyle {ht\Big (\frac
q{q^{\prime}}[\mbox {t.d.}(D)]\Big )}+\\
&&\displaystyle {\mbox {min}\Big {(}\mbox {t.d.}(\frac B{q_1}),\mbox
{t.d.}(K:D)\Big {)}}:\displaystyle {q_1\subseteq q^{\prime}\in \mbox
{ Spec}(B)\Big {\}}}+\\
&&\displaystyle {ht(\frac P{p\otimes_kB+A\otimes_kq})}$, as t.d.$(\displaystyle {\frac Ap})\leq$ t.d.$(D)\\
&&\\
&\leq&ht(p)+\mbox {max}\Big {\{}ht(q_1[\mbox {
t.d.}(A)])+\displaystyle {ht\Big (\frac {q}{q_1}[\mbox
{t.d.}(D)]\Big )}+\\
&&\displaystyle {\mbox {min}\Big {(}\mbox {t.d.}(\frac B{q_1}),\mbox
{t.d.}(K:D)\Big {)}}:\displaystyle {q_1\subseteq q\in \mbox {
Spec}(B)\Big {\}}}+\displaystyle {ht(\frac
P{p\otimes_kB+A\otimes_kq})}\\
&\leq&ht(P).
\end{array}$

\noindent Then the equality holds, as we wish to show.\\

\noindent {\bf Case 2.} $M\subseteq p^{\prime}\subset p$. By
inductive hypotheses, we get\\

$\begin{array}{lll} ht(Q)&=&ht(p^{\prime})+$max$\Big
{\{}ht(q_1[$t.d.$(A)])+ht\displaystyle {\Big (\frac
{q^{\prime}}{q_1}[\mbox {t.d.}(D)]\Big )+\mbox {min\Big {(}t.d.}(\frac
B{q_1}),\mbox {t.d.}(K:D)\Big {)}}:\\
&&q_1\subseteq q^{\prime}\in$ Spec$(B)\Big {\}}+ht(\displaystyle
{\frac Q{p^{\prime}\otimes_kB+A\otimes_kq^{\prime}})}.
\end{array}$\\

\noindent Note that $ht(p^{\prime})=ht(M)+ht(\displaystyle {\frac
{p^{\prime}}M})$ and that $\displaystyle {\frac
Q{M\otimes_kB+A\otimes_kq^{\prime}}}$ survives in
$D\otimes_kk_B(q^{\prime})$, where $k_B(q^{\prime})$ denotes the
quotient field of $\displaystyle {\frac B{q^{\prime}}}$, since
$\displaystyle {\frac Q{M\otimes_kB+A\otimes_kq^{\prime}}\cap \frac
B{q^{\prime}}}=(\overline 0)$. Then, by Proposition 2.5, we get\\

$\begin {array}{lll} ht(\displaystyle {\frac
Q{M\otimes_kB+A\otimes_kq^{\prime}}})&=&\displaystyle {ht\Big (\frac
{p^{\prime}}M\Big [\mbox {t.d.}(\frac B{q^{\prime}})}\Big ]\Big )+ht(\displaystyle
{\frac Q{p^{\prime}\otimes_kB+A\otimes_kq^{\prime}})}\\
&&\\
&&$since $\displaystyle {\frac {p^{\prime}}M\otimes_k\frac B{q^{\prime}}\cong
\frac {p^{\prime}\otimes_k B+A\otimes_kq^{\prime}}{M\otimes_kB+A\otimes_kq^{\prime}}}$, by Proposition 2.2$\\
&&\\
&=&ht(\displaystyle {\frac {p^{\prime}}M)+ht(\frac
Q{p^{\prime}\otimes_kB+A\otimes_kq^{\prime}})}\mbox { as }D\mbox {
is an AF-domain}.
\end{array}
$\\

\noindent Hence\\
$\begin{array}{lll} ht(Q)&=&ht(M)+$max$\Big
{\{}ht(q_1[$t.d.$(A)])+\displaystyle {ht\Big (\frac
{q^{\prime}}{q_1}[\mbox {t.d.}(D)]\Big )}+$min\Big {(}t.d.$\displaystyle
{(\frac
B{q_1})},$ t.d.$(K:D)\Big {)}:\\
&&q_1\subseteq q^{\prime}\in$ Spec$(B)\Big {\}}+\displaystyle {
ht(\frac Q{M\otimes_kB+A\otimes_kq^{\prime}})}.
\end{array}
$

\noindent It follows that\\

$\begin{array}{lll} ht(P)&\leq&ht(M)+$max$\Big
{\{}ht(q_1[$t.d.$(A)])+\displaystyle {ht\Big (\frac
{q^{\prime}}{q_1}}[$t.d.$(D)]\Big )+$min\Big {(}t.d.$\displaystyle
{(\frac
B{q_1})},$ t.d.$(K:D)\Big {)}:\\
&&q_1\subseteq q^{\prime}\in$ Spec$(B)\Big {\}}+\displaystyle {ht(\frac P{M\otimes_kB+A\otimes_kq^{\prime}}})\\
&&\\
&=&ht(M)+$max$\Big {\{}ht(q_1[$t.d.$(A)])+\displaystyle {ht\Big (\frac
{q^{\prime}}{q_1}}[$t.d.$(D)]\Big )+$min\Big {(}t.d.$\displaystyle
{(\frac
B{q_1})},$ t.d.$(K:D)\Big {)}:\\
&&q_1\subseteq q^{\prime}\in$ Spec$(B)\Big {\}}+\displaystyle
{ht\Big (\frac q{q^{\prime}}}[$t.d.$(D)]\Big )+\displaystyle {ht(\frac
pM)+ht(\frac P{p\otimes_kB+A\otimes_kq}}),\\
&&$ by Proposition 2.5, since $D$ is an AF-domain, and since, by Proposition 2.2,$\\
&&\\
&&\displaystyle {\frac pM\otimes_k\frac B{q^{\prime}}+\frac AM\otimes_k\frac q{q^{\prime}}\cong
\frac {p\otimes_kB+A\otimes_kq}{M\otimes_kB+A\otimes_kq^{\prime}}}\\
&&\\
&\leq&ht(p)+$max$\Big {\{}ht(q_1[$t.d.$(A)])+\displaystyle {ht\Big (\frac
q{q_1}}[$t.d.$(D)]\Big )+$min\Big {(}t.d.$\displaystyle {(\frac
B{q_1})},$
t.d.$(K:D)\Big {)}:\\
&&q_1\subseteq q\in$ Spec$(B)\Big {\}}+\displaystyle {ht(\frac P{p\otimes_kB+A\otimes_kq}})\\
&\leq&ht(P)$, as desired.$
\end{array}
$\\\\

\noindent {\bf Case 3.} $p^{\prime}\subset M$ and $q^{\prime}\neq
(0)$. Then $A_{p\prime}$ is
an AF-domain, and thus\\
$ht(Q)=ht(A\otimes_kq^{\prime})+\displaystyle {ht(\frac
Q{A\otimes_kq^{\prime}})}$, so that
$ht(P)=ht(A\otimes_kq^{\prime})+\displaystyle {ht(\frac
P{A\otimes_kq^{\prime}})}$. As $ht(\displaystyle {\frac
q{q^{\prime}})}<ht(q)$, we get, by inductive assumptions,\\

$ht(\displaystyle {\frac P{A\otimes_kq^{\prime}})}=ht(p)+\mbox {max}
\Big {\{}ht\Big (\frac {q_1}{q^{\prime}}[\mbox {t.d.}(A)]\Big )+\displaystyle
{ht\Big (\frac q{q_1}[\mbox {t.d.}(D)]\Big )+}$

\hspace {1.6cm}$\displaystyle {\mbox {min\Big {(}t.d.}(\frac
B{q_1}),\mbox {t.d.}(K:D) \Big {)}}:\displaystyle
{q^{\prime}\subseteq q_1\subseteq q\in \mbox { Spec}(B) \Big
{\}}+ht(\frac
P{p\otimes_kB+A\otimes_kq})}.$\\
\noindent Therefore\\

$ \begin{array}{lll} ht(P)&=&ht(p)+\mbox {max} \Big
{\{}ht(q^{\prime}[\mbox {t.d.}(A)])+ht\Big (\displaystyle {\frac
{q_1}{q^{\prime}}[\mbox { t.d.}(A)]\Big )}+\displaystyle {ht\Big (\frac
q{q_1}[\mbox {t.d.}(D)]\Big )+}\\
&&\mbox {min\Big {(}t.d.}(\displaystyle {\frac B{q_1}}),\mbox {
t.d.}(K:D)\Big {)}:\displaystyle {q^{\prime}\subseteq q_1\subseteq
q\in
\mbox { Spec}(B) \Big {\}}}+\\
&&\displaystyle {ht(\frac P{p\otimes_kB+A\otimes_kq})}
\end{array}$\\

$\begin{array}{lll}
&\leq&ht(p)+\mbox {max} \Big {\{}ht(q_1[\mbox
{t.d.}(A)])+\displaystyle {ht\Big {(}\frac q{q_1}[\mbox {t.d.}(D)]\Big {)}}+\\
&&\mbox {min\Big {(}t.d.}(\displaystyle {\frac B{q_1}}),\mbox {
t.d.}(K:D)\Big {)}:\displaystyle {q_1\subseteq q\in
\mbox { Spec}(B) \Big {\}}+ht(\frac P{p\otimes_kB+A\otimes_kq})}\\
&\leq&ht(P)$, and the equality holds.$
\end{array}
$\\\\

\noindent {\bf Case 4.} $(0)\neq p^{\prime}\subset M$ and
$q^{\prime}=(0)$. Then $ht(Q)=\displaystyle
{ht(p^{\prime}\otimes_kB)+ht(\frac Q{p^{\prime}\otimes_kB}})$, via Proposition 2.5. It
follows that $ht(P)=\displaystyle {ht(p^{\prime})+ht(\frac
P{p^{\prime}\otimes_kB}})$, as $A_{p{\prime}}$ is an AF-domain. Since $p^{\prime}\subset M$, there exists a
unique prime ideal $P^{\prime}$ of $T$ such that $P^{\prime}\cap
A=p^{\prime}$. Then $P^{\prime}\neq (0)$, so that
dim$(\displaystyle {\frac T{P^{\prime}})}<$ dim$(T)$, and
$$\begin{array}{lll}
\displaystyle {\frac A{p^{\prime}}}&\longrightarrow& D\\
&&\\
\downarrow&&\downarrow\\
&&\\
\displaystyle {\frac T{P^{\prime}}}&\longrightarrow&\displaystyle
{K\cong \frac {T/P^{\prime}}{M/P^{\prime}}}
\end{array}
$$ is a pullback diagram.
Moreover, by Proposition 2.6, $\displaystyle {\frac T{P^{\prime}}}$ is an AF-domain which is catenarian,
therefore, by inductive
assumptions, we get
$$ht(\displaystyle {\frac P{p^{\prime}\otimes_kB})=ht(\frac
p{p^{\prime}})+\mbox {max}\Big {\{}ht\Big (q_1\Big [\mbox {t.d.}(\frac
A{p^{\prime}})\Big ]\Big )}+\displaystyle {ht\Big (\frac q{q_1}[\mbox
{t.d.}(D)]\Big )}+$$ $$\mbox {min\Big {(}t.d.}(\displaystyle
{\frac B{q_1}),\mbox {t.d.}(K:D)\Big {)}}:\displaystyle
{q_1\subseteq q\in \mbox { Spec}(B)\Big {\}}+ht\Big {(}\frac
{P/(p^{\prime}\otimes_kB)}{p/p^{\prime}\otimes_kB+A/p^{\prime}\otimes_kq}\Big
{)}}$$
$$=ht(\frac p{p^{\prime}})+\mbox {max}\Big
{\{}ht\Big (q_1\Big [\mbox {t.d.}(\frac A{p^{\prime}})\Big ]\Big )+\displaystyle
{ht\Big (\frac q{q_1}[\mbox {t.d.}(D)]\Big )}+$$
$$\mbox {min\Big {(}t.d.}(\frac B{q_1}),\mbox {t.d.}(K:D)\Big
{)}:\displaystyle {q_1\subseteq q\in \mbox { Spec}(B)\Big
{\}}+ht(\frac {P}{p\otimes_kB+A\otimes_kq})}.$$ Hence
$$\begin{array}{lll}
ht(P)&=&ht(p^{\prime})+ht(\displaystyle {\frac p{p^{\prime}})+\mbox
{max}\Big {\{}ht\Big (q_1\Big [\mbox { t.d.}(\frac
A{p^{\prime}})\Big ]\Big )}+\displaystyle {ht\Big (\frac q{q_1}[\mbox {t.d.}(D)]\Big )}+\\
&&\mbox {min}\Big {(}\mbox {t.d.}(\displaystyle {\frac
B{q_1}}),\mbox {t.d.}(K:D)\Big {)}:q_1\subseteq q\in \mbox {
Spec}(B) \Big {\}}+ht(\displaystyle {\frac {P}{p\otimes_kB+A\otimes_kq})}\\
&&\\
&\leq&ht(p)+\mbox {max}\Big {\{}ht(q_1[\mbox { t.d.}(
A)])+\displaystyle {ht\Big (\frac q{q_1}[\mbox {t.d.}(D)]\Big
)+}\\
&&\mbox {min\Big (t.d.}(\displaystyle {\frac B{q_1}}),\mbox
{t.d.}(K:D)\Big ):\displaystyle {q_1\subseteq q\in \mbox {
Spec}(B)\Big {\}}
+ht(\frac {P}{p\otimes_kB+A\otimes_kq})}\\
&\leq&h(P).
\end{array}$$ Then the equality holds.\\

\noindent {\bf Case 5.} $p^{\prime}=(0)$ and $q^{\prime}=(0)$. Then, by [19, Theorem 3.1], $ht(Q)\leq$
t.d.$(B)$, and thus, as $Q_{r+1}\cap A=p$ and $Q_{r+1}\cap B=q$, we get

$$ht(P)=ht(Q_{r+1})+\displaystyle {ht(\frac P{Q_{r+1}})\leq 1+\mbox
{t.d.}(B)+ht(\frac P{p\otimes_kB+A\otimes_kq})}.$$ Suppose that
$1+$t.d.$(B)\leq ht(p[$t.d.$(B)])$. Therefore,\\

$\begin{array}{lll} ht(P)&\leq& ht(p[$t.d.$(B)])+ht(\displaystyle
{\frac P{p\otimes_kB+A\otimes_kq}})\\
&&\\
&\leq&ht(p[$t.d.$(B)])+ht\Big (q\Big [$t.d.$(\displaystyle {\frac
Ap})\Big ]\Big )+ht(\displaystyle {\frac P{p\otimes_kB+A\otimes_kq}})\\
&&\\
&=& ht(p\otimes_kB)+ht(\displaystyle {\frac
{p\otimes_kB+A\otimes_kq}{p\otimes_kB})}+ht(\displaystyle {\frac
P{p\otimes_kB+A\otimes_kq}})\\
&\leq&ht(P).
\end{array}$\\

\noindent Hence
$ht(p[$t.d.$(B)])=ht(p[$t.d.$(B)])+ht\Big (q\Big [$t.d.$(\displaystyle {\frac
Ap)\Big ]\Big )}$, so that $q=(0)$ which leads to a contradiction, as $ht(q)>0$. It
follows that, by Proposition 2.3,\\

\hspace{1cm} $1+$t.d.$(B)>ht(p[$t.d.$(B)])=ht(p)+$min\Big
(t.d.$(B),$
t.d.$(K:D)\Big )$, so that\\

\hspace{2.5cm} t.d.$(B)>$ t.d.$(K:D)$, as $ht(p)\geq 1$. On the
other hand,

$\begin{array}{lll} ht(P)&=&ht(Q_{r+1})+ht(\displaystyle {\frac
P{Q_{r+1}}})\\
&=&\lambda \Big (\Big ((0),(0)\Big ),Q_{r+1}\Big {)}+ht(\displaystyle {\frac
P{Q_{r+1}}})\\
&\leq&$t.d.$(A)-$t.d.$(\displaystyle {\frac
Ap})+ht\Big (q\Big [$t.d.$(\displaystyle {\frac Ap)\Big ]\Big )+ht(\frac
{Q_{r+1}}{p\otimes_kB+A\otimes_kq})+ht(\frac P{Q_{r+1}})}\\
&&\mbox {(cf. Proposition 2.7)}\\
&&\\
&\leq& ht(p[$t.d.$(K:D)])+ht\Big (q\Big [$t.d.$(\displaystyle {\frac
Ap)\Big ]\Big )+ht(\frac P{p\otimes_kB+A\otimes_kq})},\\
&&\mbox { as }A[\mbox {t.d.}(K:D)]\mbox { is an AF-domain}\\
&\leq& ht(p[$t.d.$(B)])+ht\Big (q[$t.d.$(\displaystyle {\frac
Ap)]\Big )+ht(\frac P{p\otimes_kB+A\otimes_kq})}
\end{array}$\\

$\begin{array}{lll}
&\leq&ht(p)+\mbox {max}\Big {\{}ht(q_1[\mbox {
t.d.}(A)])+\displaystyle {ht\Big (\frac q{q_1}[\mbox
{t.d.}(D)]\Big )}+\\
&&\mbox {min\Big {(}t.d.}(\displaystyle {\frac B{q_1}}),\mbox
{t.d.}(K:D)\Big {)}:\displaystyle {q_1\subseteq q\in \mbox {
Spec}(B)\Big
{\}}+ht(\frac P{p\otimes_kB+A\otimes_kq})},\\
&&\mbox { by }(*)\\
&\leq& ht(P).
\end{array}
$\\

\noindent Then the equality holds.\\

\noindent {\bf \underline {Step 2}.} $B$ is an arbitrary $k$-algebra.\\

\noindent Let $P_0\subset P_1\subset...\subset P_h=P$ be a chain of
prime ideals of $A\otimes_kB$ such that $h=ht(P)$. Let $q_0:=P_0\cap
B$. Then
$$\displaystyle {\frac {P_0}{A\otimes_kq_0}\subset \frac
{P_0}{A\otimes_kq_0}\subset \frac {P_1}{A\otimes_kq_0}\subset
...\subset \frac {P_h}{A\otimes_kq_0}=\frac {P}{A\otimes_kq_0}}$$ is
a chain of prime ideals of $A\otimes_k\displaystyle {\frac B{q_0}}$
and $h=ht(P)=ht(\displaystyle {\frac P{A\otimes_kq_0})}$. By Step 1,
$$\begin{array}{lll}
ht(P)&=&ht(\displaystyle {\frac P{A\otimes_kq_0})=ht(p)+\mbox
{max}\Big {\{}ht\Big (\frac {q_1}{q_0}[\mbox { t.d.}(A)]\Big )}+\displaystyle
{ht\Big (\frac
q{q_1}[\mbox {t.d.}(D)]\Big )}+\\
&&\\
&&\displaystyle {\mbox {min}\Big {(}\mbox {t.d.}(\frac B{q_1}),\mbox
{ t.d.}(K:D)\Big {)}}:\displaystyle {q_0\subseteq q_1\subseteq q\in \mbox {
Spec}(B)\Big {\}}+}\\
&&\\
&&\displaystyle {ht\Big (\frac {P/(A\otimes_kq_0)}{p\otimes_k(B/q_0)+A\otimes_k(q/q_0)}\Big )}
\end{array}$$
$$\begin{array}{lll}
&\leq&ht(p)+\mbox {max}\Big {\{}ht(q_1[\mbox {
t.d.}(A)])+\displaystyle {ht\Big (\frac
q{q_1}[\mbox {t.d.}(D)]\Big )}+\\
&&\\
&&\displaystyle {\mbox {min}\Big {(}\mbox {t.d.}(\frac B{q_1}),\mbox
{ t.d.}(K:D)\Big {)}}:\displaystyle {q_1\subseteq q\in \mbox {
Spec}(B)\Big {\}}+}\displaystyle {ht(\frac {P}{p\otimes_k B+A\otimes_k q})}\\
&\leq&ht(P),\mbox { then the equality holds establishing the desired
formula}.
\end{array}$$

3) First, observe that, by [1, Lemma 2.1], for each $p\in$ Spec$(A)$ such that $M\not\subset p$, there
exists a unique $Q\in$ Spec$(T)$ such that $Q\cap A=p$, and $Q$ satisfies $A_p=T_Q$. Then\\
$d=$ max$\Big \{ht(p):p\in$ Spec$(A)$ with $M\not\subset p\Big \}$. Now, by (1), we have

max$\Big \{ht(P):P\in$ Spec$(A\otimes_kB)$ and $M\not\subset p:=P\cap A\Big \}=$

max$\Big \{ht(p)+ht(q[$t.d.$(A)])+ht(\displaystyle {\frac P{p\otimes_kB+A\otimes_kq}}):P\in$ Spec$(A\otimes_kB)$
with $p=P\cap A$,\\

\hspace {2cm} $q=P\cap B$ such that $M\not\subset p\Big \}=$

max$\Big \{ht(p)+ht(q[$t.d.$(A)])+\displaystyle {\mbox {min\Big (t.d.(}\frac Ap),\mbox { t.d.}(\frac Bq)}\Big ):
p\in$ Spec$(A)$ and $q\in$ Spec$(B)$\\

\hspace {2cm}such that $M\not\subset p\Big \}$ (cf. [20, Proposition 2.3]) $=$

max$\Big \{ht(q[$t.d.$(A)])+\displaystyle {\mbox {min\Big (t.d.(}A),ht(p)+\mbox {t.d.}(\frac Bq)}\Big ):p\in$
Spec$(A)$
with $M\not\subset p$ and\\

\hspace {2cm} $q\in$ Spec$(B)\Big \}$ (as $A_p$ is an AF-domain) $=$

max$\Big \{ht(q[$t.d.$(A)])+\displaystyle {\mbox {min\Big (t.d.(}A),d+\mbox {t.d.}(\frac Bq)}\Big ):q\in$
Spec$(B)\Big \}=D\Big ($t.d.$(A),d,B\Big )$ $(**)$. On the other hand, let $M\subseteq p$. As done in (2), we may
assume that $(T,M)$ is a quasilocal domain, and thus $M$ is a divided prime ideal of $A$. First, note that\\

$\begin{array}{lll}
ht(p)+
$min\Big (t.d.$(\displaystyle {\frac Ap),\mbox { t.d.}(\frac Bq)\Big )}&=&\displaystyle {ht(M)+ht(\frac pM)+
\mbox {min}\Big (\mbox {t.d.}
(\frac Ap),\mbox { t.d.}(\frac Bq)\Big )}\\
&=&\displaystyle {ht(M)+\mbox {min\Big (t.d.}(D),ht(\frac pM)+\mbox {t.d.}(\frac Bq)\Big )}$, as $D$ is$\\
&&$an AF-domain$.
\end{array}$\\

\noindent Hence, by (2),\\

max$\Big \{ht(P):P\in$ Spec$(A\otimes_kB)$ and $M\subseteq p:=P\cap
A\Big \}=$\\

$\mbox {max}\Big \{ht(p)+\mbox {max}\Big \{ht(q_1[\mbox
{t.d.}(A)])+\displaystyle {ht\Big (\frac q{q_1}[\mbox {t.d.}(D)]\Big
)}+\displaystyle {\mbox {min}\Big {(}\mbox {t.d.}(\frac
B{q_1}),\mbox { t.d.}(K:D)\Big {)}}:$ $$\displaystyle {q_1\subseteq
q\in \mbox { Spec}(B)\Big {\}}+ht(\frac
P{p\otimes_kB+A\otimes_kq})}:P\in\mbox { Spec}(A\otimes_kB)\mbox {
with }p=P\cap A \mbox { and }$$ $$ q=P\cap B\mbox { such that }
M\subseteq p\Big \}=$$

$\mbox {max}\Big \{ht(p)+ht(q_1[\mbox {t.d.}(A)])+\displaystyle
{ht\Big (\frac q{q_1}[\mbox {t.d.}(D)]\Big )}+\displaystyle {\mbox
{min}\Big {(}\mbox {t.d.}(\frac B{q_1}),\mbox { t.d.}(K:D)\Big
{)}}+$ $$\mbox {min}\Big (\mbox {t.d.}(\frac Ap),\mbox { t.d.}(\frac
Bq)\Big ):p\in\mbox { Spec}(A)\mbox { and }q_1\subseteq q\in \mbox {
Spec}(B)\mbox { such that }M\subseteq p \Big \}=$$

$\mbox {max}\Big \{ht(M)+ht(q_1[\mbox {t.d.}(A)])+\displaystyle
{ht\Big (\frac q{q_1}[\mbox {t.d.}(D)]\Big )}+\displaystyle {\mbox
{min}\Big {(}\mbox {t.d.}(\frac B{q_1}),\mbox { t.d.}(K:D)\Big
{)}}+$ $$\mbox {min}\Big (\mbox {t.d.}(D),ht(\frac pM)+\mbox {t.d.}(\frac Bq)\Big ):
p\in\mbox { Spec}(A)\mbox { and }q_1\subseteq
q\in \mbox { Spec}(B)\mbox { such that }M\subseteq p \Big \}=$$

$\mbox {max}\Big \{ht(M)+ht(q_1[\mbox {t.d.}(A)])+\displaystyle
{ht\Big (\frac q{q_1}[\mbox {t.d.}(D)]\Big )}+\displaystyle {\mbox
{min}\Big {(}\mbox {t.d.}(\frac B{q_1}),\mbox { t.d.}(K:D)\Big
{)}}+$ $$\mbox {min}\Big (\mbox {t.d.}(D),\mbox { dim}(D)+\mbox {t.d.}(\frac Bq)\Big ):
q_1\subseteq q\in \mbox { Spec}(B)\Big \}\mbox
{  }(***).$$

\noindent Consequently, combining $(**)$ and $(***)$, we get easily the
desired formula for dim$(A\otimes_kB)$.

4) Let $p\in$ Spec$(A)$ and $q\in$ Spec$(B)$. Let $Q$ be a minimal prime ideal of $p\otimes_kB+A\otimes_kq$.
Assume that $M\not\subset p$. Then, by (1),
$ht(Q)=ht(p)+ht(q[$t.d.$(A)])$. Therefore, $$ht(p\otimes_kB+A\otimes_kq)=ht(p)+ht(q[\mbox {t.d.}(A)])$$ and for any
$P\in$ Spec$(A\otimes_kB)$ with $p=P\cap A$ and $q=P\cap B$, via (1),
$$ht(P)=ht(p\otimes_kB+A\otimes_kq)+ht(\displaystyle {\frac P{p\otimes_kB+A\otimes_kq}}).$$
Now, let $M\subseteq p$. Then, by (2), $$ht(Q)=ht(p)+\mbox {max}\Big
{\{}ht(q_1[\mbox { t.d.}(A)])+\displaystyle {ht\Big (\frac
q{q_1}[\mbox {t.d.}(D)]\Big )}+\displaystyle {\mbox {min}\Big {(}\mbox {t.d.}(\frac B{q_1}),\mbox
{ t.d.}(K:D)\Big {)}}:$$ $$\displaystyle {q_1\subseteq q\in \mbox { Spec}(B)\Big {\}}}, \mbox { so that}$$
$$ht(p\otimes_kB+A\otimes_kq)=ht(p)+\mbox {max}\Big
{\{}ht(q_1[\mbox { t.d.}(A)])+\displaystyle {ht\Big (\frac
q{q_1}[\mbox {t.d.}(D)]\Big )}+$$ $$\displaystyle {\mbox {min}\Big {(}\mbox {t.d.}(\frac B{q_1}),
\mbox
{ t.d.}(K:D)\Big {)}}:\displaystyle {q_1\subseteq q\in \mbox {
Spec}(B)\Big {\}}}.$$ Thus, via (2), $ht(P)=ht(p\otimes_kB+A\otimes_kq)+ht(\displaystyle {\frac P{p\otimes_kB+
A\otimes_kq}})$
for each $P\in$ Spec$(A\otimes_kB)$ such that $p=P\cap A$ and $q=P\cap B$.
Therefore $A$ is a GAF-domain, completing the proof.

\hspace {13.3cm}$\Box$\\

We get the following interesting consequence of Theorem 2.8.\\

\noindent {\bf Corollary 2.9.} {\it Let $T$ be a $k$-algebra which
is a domain and $M$ a maximal ideal of $T$. Let $K=\displaystyle
{\frac TM}$ and $D$ be a subring of $K$.
Let $\varphi:T\longrightarrow K$ be the canonical surjective
homomorphism and $A:=\varphi^{-1}(D)$. Assume that $T$ and $D$ are
AF-domains and $ht(M)\leq 2$. Then, the assertions (1), (2), (3) and (4) of Theorem 2.8 hold for $A$ and
any $k$-algebra $B$.}\\

\noindent {\bf Proof.} It is direct from Theorem 2.8 since any
two-dimensional domain is catenarian. $\Box$\\

The following result discusses a case where it is possible to drop the catenarity assumption of $T_M$ in Theorem 2.8.\\

\noindent {\bf Proposition 2.10.} {\it Let $T$ be a $k$-algebra which is a domain and
$M$ a maximal ideal of $T$. Let $K=\displaystyle {\frac TM}$ and $D$ be a subring of $K$. Let
$\varphi:T\longrightarrow K$ be the canonical surjective
homomorphism and $A:=\varphi^{-1}(D)$. Assume that $T$ and $D$ are
AF-domains and t.d.$(K:D)\leq 2$. Then, the assertions (1), (2), (3) and (4) of Theorem 2.8 hold for $A$ and
any $k$-algebra $B$.}\\

\noindent {\bf Proof.} The proof runs similar to that of Theorem 2.8.
We need only to discuss the case where $(0)\neq p^{\prime}\subset M$ and
$q^{\prime}=(0)$ (Case 4 of Step 1 of the above proof) where the catenarity property of $T$ is required.
So, let $B$ be an integral domain and let $P\in$ Spec$(A\otimes_kB)$, $p=P\cap A$ and $q=P\cap B$. Let
$r:=$ max$\{m:Q_m\cap A\subset p$ or $Q_m\cap B\subset q\}$. Let
$Q=Q_r$, $p^{\prime}=Q_r\cap A$ and $q^{\prime}=Q_r\cap B$, and
suppose that $(0)\neq p^{\prime}\subset M$ and
$q^{\prime}=(0)$. Hence\\

$\begin{array}{lll} ht(P)&=&ht(Q_{r+1})+ht(\displaystyle {\frac PQ_{r+1}})\\
&=&\lambda \Big {(}(.,(0)),Q_{r+1}\Big {)}+ht(\displaystyle {\frac PQ_{r+1}})\\
&\leq& $t.d.$(A)-$t.d.$(\displaystyle {\frac Ap)+ht\Big (q[\mbox
{t.d.}(\frac Ap)]\Big )+ht(\frac {Q_{r+1}}{p\otimes_kB+A\otimes_kq})}+ht(\displaystyle {\frac PQ_{r+1}}),\\
&&$ by Proposition 2.7$\\
&\leq&ht(p[$t.d.$(K:D)])+\displaystyle {ht\Big (q[\mbox {t.d.}(\frac
Ap)]\Big )+ht(\frac P{p\otimes_kB+A\otimes_kq})}$, as $A[$t.d.$(K:D)]$ is$\\
&&$an AF-domain, by Proposition 2.3$.
\end{array}
$\\

\noindent Now, if t.d.$(B)\geq 2$, then we get
$$
\begin{array}{lll}ht(P)&\leq& ht(p[\mbox {t.d.}(B)])+\displaystyle {ht\Big (q[\mbox
{t.d.}(\frac Ap)]\Big )+ht(\frac P{p\otimes_kB+A\otimes_kq})}\\
&\leq& ht(p)+\mbox {max}\Big {\{}ht(q_1[\mbox {
t.d.}(A)])+\displaystyle {ht\Big (\frac q{q_1}[\mbox {t.d.}(D)]\Big )+\mbox
{min\Big {(}t.d.}(\frac B{q_1}),\mbox {t.d.}(K:D)\Big {)}}:\\
&&\displaystyle {q_1\subseteq q\in \mbox { Spec}(B)\Big
{\}}+ht(\frac P{p\otimes_kB+A\otimes_kq})}\\
&\leq& ht(P), \end{array}$$ so that the desired equality holds.
Assume that
t.d.$(B)\leq 1$. Then $B$ is an AF-domain. Hence\\

$\begin{array}{lll} ht(P)&=&ht(p[$t.d.$(B)])+ht\Big (q[\mbox
{t.d.}\displaystyle {(\frac
Ap)]\Big )+ht(\frac P{p\otimes_kB+A\otimes_kq}})$, by Proposition 2.5$\\
&\leq&ht(p)+\mbox {max}\Big {\{}ht(q_1[\mbox {
t.d.}(A)])+\displaystyle {ht\Big (\frac q{q_1}[\mbox {t.d.}(D)]\Big )+\mbox
{min\Big {(}t.d.}(\frac B{q_1}),\mbox {t.d.}(K:D)\Big {)}}:\\
&&\displaystyle {q_1\subseteq q\in \mbox { Spec}(B)\Big
{\}}+ht(\frac P{p\otimes_kB+A\otimes_kq})}\\
&\leq& ht(P).
\end{array}
$\\
Then the equality holds, as contended. $\Box$\\

Next, for each positive integer $n$, we provide an example of a finite dimensional valuation domain $(V,M)$, thus universally
catenarian (cf. [18]), of Krull dimension $n$ such that
$V$ is an AF-domain. Then, given any subring $D$ of $\displaystyle {K:=\frac VM}$ such that $D$ is an AF-domain (in
particular any subfield of $K$), Theorem 2.8 allows
us to compute dim$(A\otimes_kB)$ for any $k$-algebra $B$, where $A$ is the pullback $\varphi^{-1}(D)$.\\

\noindent {\bf Example 2.11.} Let
$V_1=k(X_1,X_2,...,X_n)[Y]_{(Y)}=k(X_1,X_2,...X_n)+M_1$ be a
rank-one discrete valuation domain with $M_1:=YV_1$. Let
$$V_2=k(X_1,X_2,...,X_{n-1})[X_n]_{(X_n)}+M_1=k(X_1,X_2,...,X_{n-1})+M_2$$
with $M_2:=M_1+X_nk(X_1,...,X_{n-1})[X_n]_{(X_n)}$. It is well known that
$V_2$ is a valuation domain of (Krull) dimension $2$ [11, Exercise 13 (2), page 203]. Also, by [12],
$V_2$ is an AF-domain. Let
$$V_3=k(X_1,X_2,...,X_{n-2})[X_{n-1}]_{(X_{n-1})}+M_2=k(X_1,X_2,...,X_{n-2})+M_3$$
with $M_3:=M_2+X_{n-1}k(X_1,...,X_{n-2})[X_{n-1}]_{(X_{n-1})}$. Then $V_3$ is a valuation domain
of dimension $3$ which is an
AF-domain.  We may iterate this process to construct a valuation
domain of the form $V=K+M$ which is an AF-domain of dimension $r$
for each positive integer $r$. $\Box$\\

\noindent {\bf Acknowledgement.} The author would like to thank the referee for his/her helpful suggestions.\\

\noindent {\bf References}
\begin{list}{}{\topsep=2mm \parsep=0mm \itemsep=0mm \leftmargin=7.1mm}
{\small

\item [{[1]}] D. F. Anderson, A. Bouvier, D. E. Dobbs, M. Fontana, and S. Kabbaj,
{\it On Jaffard domains}, Exposition. Math. 6 (2) (1988) 145-175.
\item [{[2]}]  S. Bouchiba, {\it On Krull dimension of tensor products of algebras arising from AF-domains},
J. Pure Appl. Algebra 203 (2005) 237-251.
\item [{[3]}]  S. Bouchiba, {\it Chains of prime ideals in tensor products of algebras},
J. Pure Appl. Algebra 209 (2007) 621-630.
\item [{[4]}]  S. Bouchiba, {\it AF-domains and Jaffard's special
chain theorem}, The Arabian Journal for Science and Engineering,
Vol. 30, n: 2A, (2005), 1-7.
\item [{[5]}]  S. Bouchiba, {\it AF-domains and their
generalizations}, Comm. Algebra, to appear.
\item [{[6]}]  S. Bouchiba, F. Girolami, and S. Kabbaj, {\it The dimension of tensor
products of k-algebras arising from pullbacks}, J. Pure Appl.
Algebra 137 (1999), 125-138.
\item [{[7]}]  S. Bouchiba and S. Kabbaj, {\it Tensor products of Cohen-Macaulay rings: Solution to a problem of
Grothendieck}, J. Algebra 252 (2002) 65-73.
\item [{[8]}] S. Bouchiba, D. Dobbs, and S. Kabbaj, {\it On the prime ideal structure of
tensor products of algebras}, J. Pure Appl. Algebra 176 (2002), 89-112.
\item [{[9]}] A. Bouvier, D. Dobbs
and M. Fontana, {\it Universally catenarian integral domains}, Adv. Math. 72 (1988), 211-238.
\item [{[10]}] J.W. Brewer, P.R. Montgomery, E.A. Rutter, W.J.
Heinzer, {\it Krull dimension of polynomial rings}, Lecture Notes in
Math., vol. 311, Springer-Verlag, Berlin, NY, 1972, 26-45.
\item [{[11]}]  R. Gilmer, {\it Multiplicative ideal theory}, Marcel Dekker, New
York, 1972.
\item [{[12]}]  F. Girolami, {\it AF-Rings and locally Jaffard
Rings}, Lecture Notes Pure Appl. Math., Vol. 153. New York: Dekker,
1994, 151-161.
\item [{[13]}]  A. Grothendieck, {\it El\'ements de g\'eom\'etrie
alg\'ebrique}, Vol. IV,  Institut des Hautes Etudes Sci. Publ. Math.
No. 24, Bures-sur-Yvette, 1965.
\item [{[14]}] P. Jaffard, {\it Th\'eorie de la dimension dans les anneaux de polyn\^omes}, M\'em. Sc. Math.
146, Gauthier-Villars, Paris, 1960.
\item [{[15]}] I. Kaplansky, {\it Commutative rings}, Chicago University Press, 1974.
\item [{[16]}] H. Matsumura, {\it Commutative ring theory}, Cambridge University
Press, 1986.
\item [{[17]}] M. Nagata, {\it Local rings}, Interscience, New York,
1962.
\item [{[18]}] M. Nagata, {\it Finitely generated rings over a
valuation domain}, J. Math. Kyoto Univ. 5 (1966), 163-169.
\item [{[19]}]  R.Y. Sharp, {\it The dimension of the tensor product of two field
extensions}, Bull. London Math. Soc. 9 (1977), 42-48.
\item [{[20]}]  A.R. Wadsworth, {\it The Krull dimension of tensor products of
commutative algebras over a field}, J. London Math. Soc. 19 (1979),
391-401.
\item [{[21]}]  O. Zariski and P. Samuel, {\it Commutative algebra Vol. II}, Van
Nostrand, Princeton, 1960.}

\end{list}
\end{document}